%% file: altdistn.tex
\pdfoutput=1

\documentclass[reprint,amsmath,amssymb,floatfix,longbibliography,aip,letterpaper,cha]{revtex4-1}

\usepackage{graphics,graphicx,float}

\makeatletter
\preprintsty@sw{}{%
 \usepackage{titlesec}
 \titleformat{\section}{\normalfont\footnotesize\sffamily\bfseries\uppercase}%
 	{\thesection}{1em}{}%
 \titleformat{\subsection}{\normalfont\small\sffamily\bfseries}%
 	{\thesubsection}{1em}{}%
 \titleformat{\subsubsection}{\normalfont\small\sffamily\slshape}{\thesubsubsection}{1em}{}%
 }%
\makeatother
  
\def\citeyear{\citep}
\def\autocite{\citep}

\usepackage[pdftex,
            pdfauthor={Steven A. Frank},
            pdftitle={Invariance defines scaling relations and probability patterns},
            pdfsubject={Invariance defines scaling relations and probability patterns}]{hyperref} 

\input altdistn-def

\usepackage[numbib,notindex,nottoc,notlof,notlot]{tocbibind}
\settocbibname{References}

\newif\ifshowgit \showgittrue		
\newif\ifshowtime \showtimetrue		
\ifshowgit
\usepackage[calc,timezonesep=, showseconds=false]{datetime2} 	
\DTMsavenow{now}
\DTMtozulu{now}{zdate}			
\usepackage[grumpy,mark,dirty=-*,raisemark=0.02\paperheight,local]{gitinfo2}
\renewcommand{\gitMarkPref}{git}
\renewcommand{\gitMark}{\gitBranch\,@\,\gitRel-\gitRoff-\gitAbbrevHash\gitDirty{}-%
	\gitAuthorDate\ \ifshowtime(\DTMuse{zdate})\fi\ \textbullet{}\ safrank}
\fi

\usepackage{titlesec}

\renewcommand{\thesection}{\arabic{section}}
\renewcommand{\thesubsection}{\thesection.\arabic{subsection}}
\renewcommand{\thesubsubsection}{\thesubsection.\arabic{subsubsection}}

\begin{document}

\title{Common probability patterns arise from simple invariances\\ \phantom{x}}

\author{Steven A.\ Frank}
\affiliation{Department of Ecology and Evolutionary Biology, University of California, Irvine, CA 92697--2525  USA}

\begin{abstract}

Shift and stretch invariance lead to the exponential-Boltzmann probability distribution. Rotational invariance generates the Gaussian distribution. Particular scaling relations transform the canonical exponential and Gaussian patterns into the variety of commonly observed patterns. The scaling relations themselves arise from the fundamental invariances of shift, stretch, and rotation, plus a few additional invariances. Prior work described the three fundamental invariances as a consequence of the equilibrium canonical ensemble of statistical mechanics or the Jaynesian maximization of information entropy. By contrast, I emphasize the primacy and sufficiency of invariance alone to explain the commonly observed patterns. Primary invariance naturally creates the array of commonly observed scaling relations and associated probability patterns, whereas the classical approaches derived from statistical mechanics or information theory require special assumptions to derive commonly observed scales.  
\phantom{\footnote{web: \href{http://stevefrank.org}{http://stevefrank.org}}}

\bigskip


\end{abstract}

\maketitle

{\renewcommand{\tocname}{}\small\hbox{\null}\vskip-66pt\tableofcontents}\newpage


\begin{quote}
\small{
\baselineskip=13pt

It is increasingly clear that the symmetry [invariance] group of nature is the deepest thing that we understand about nature today. I would like to suggest something here that I am not really certain about but which is at least a possibility: that specifying the symmetry group of nature may be all we need to say about the physical world, beyond the principles of quantum mechanics.

The paradigm for symmetries of nature is of course the group symmetries of space and time. These are symmetries that tell you that the laws of nature don't care about how you orient your laboratory, or where you locate your laboratory, or how you set your clocks or how fast your laboratory is moving (\textcite[p.~73]{weinberg99towards}).

\centerline{\rule{80pt}{1.0pt}}

For the description of processes taking place in nature, one must have a \emph{system of reference} (\textcite[p.~1]{landau80mechanics}).

}
\end{quote}

\section{Introduction}

I argue that three simple invariances dominate much of observed pattern. First, probability patterns arise from invariance to a shift in scaled measurements. Second, the scaling of measurements satisfies  invariance to uniform stretch. Third, commonly observed scales are often invariant to rotation.

Feynman \autocite{feynman98statistical} described the shift invariant form of probability patterns as
\begin{equation}\label{eq:ratioInv}
  \frac{q\lrp{\E}}{q\lrp{\Ep}} = \frac{q\lrp{\E+a}}{q\lrp{\Ep+a}},
\end{equation}
in which $q\lrp{\E}$ is the probability associated with a measurement, $\E$. Here, the ratio of probabilities for two different measurements, $\E$ and $\Ep$, is invariant to a shift by $a$. Feynman derived this invariant ratio as a consequence of Boltzmann's equilibrium distribution of energy levels, $\E$, that follows from statistical mechanics
\begin{equation}\label{eq:boltz}
  q\lrp{\E} = \Gl e^{-\Gl\E}.
\end{equation}
Here, $\Gl=1/\angb{\E}$ is the inverse of the average measurement. 

Feynman presented the second equation as primary, arising as the equilibrium from the underlying dynamics of particles and the consequent distribution of energy, $\E$. He then mentioned in a footnote that the first equation of shift invariance follows as a property of equilibrium. However, one could take the first equation of shift invariance as primary. The second equation for the form of the probability distribution then follows as a consequence of shift invariance.

What is primary in the relation between these two equations: equilibrium statistical mechanics or shift invariance? The perspective of statistical mechanics, with \Eq{boltz} as the primary equilibrium outcome, dominates treatises of physics.

Jaynes \autocite{jaynes57information,jaynes57informationII} questioned whether statistical mechanics is sufficient to explain why patterns of nature often follow the form of \Eq{boltz}. Jaynes emphasized that the same probability pattern often arises in situations for which physical theories of particle dynamics make little sense. In Jaynes' view, if most patterns in economics, biology, and other disciplines follow the same distributional form, then that form must arise from principles that transcend the original physical interpretations of particles, energy, and statistical mechanics \autocite{jaynes03probability}. 

Jaynes argued that probability patterns derive from the inevitable tendency for systems to lose information. By that view, the equilibrium form expresses minimum information, or maximum entropy, subject to whatever constraints may act in particular situations. In maximum entropy, the shift invariance of the equilibrium distribution is a consequence of the maximum loss of information under the constraint that total probability is conserved. 

Here, I take the view that shift invariance is primary. My argument is that shift invariance and the conservation of total probability lead to the exponential-Boltzmann form of probability distributions, without the need to invoke Boltzmann's equilibrium statistical mechanics or Jaynes' maximization of entropy. Those secondary special cases of Boltzmann and Jaynes follow from primary shift invariance and the conservation of probability. The first part of this article develops the primacy of shift invariance. 

Once one adopts the primacy of shift invariance, one is faced with the interpretation of the measurement scale, $\E$. We must abandon energy, because we have discarded the primacy of statistical mechanics, and we must abandon Jaynes' information, because we have assumed that we have only general invariances as our basis. 

We can of course end up with notions of energy and information that derive from underlying invariance. But that leaves open the problem of how to define the canonical scale, $\E$, that sets the frame of reference for measurement. 

We must replace the scaling relation $\E$ in the above equations by something that derives from deeper generality: the invariances that define the commonly observed scaling relations. 

In essence, we start with an underlying scale for observation, $z$. We then ask what transformed scale, $z\mapsto\trz\equiv\E$, achieves the requisite shift invariance of probability pattern, arising from the invariance of total probability. It must be that shift transformations, $\trz\mapsto a+ \trz$, leave the probability pattern invariant, apart from a constant of proportionality. 

Next, we note that a stretch of the scale, $\trz\mapsto b\trz$, also leaves the probability pattern unchanged, because the inverse of the average value in \Eq{boltz} becomes $\Gl=1/b\angb{\trz}$, which cancels the stretch in the term $\Gl\E=\Gl\trz$. Thus, the scale $\trz$ has the property that the associated probability pattern is invariant to the affine transformation of shift and stretch, $\trz\mapsto a+b\trz$. That affine invariance generates the symmetry group of scaling relations that determine the commonly observed probability patterns \autocite{frank10measurement,frank11a-simple,frank14how-to-read}.

The final part of this article develops rotational invariance of conserved partitions. For example, the Pythagorean partition $\trz=x^2(s) + y^2(s)$ splits the scaled measurement into components that add invariantly to $\trz$ for any value of $s$. The invariant quantity defines a circle in the $xy$ plane with a conserved radius $\Rz=\sqrt{\trz}$ that is invariant to rotation around the circle, circumscribing a conserved area $\pi\Rz^2=\pi\trz$. Rotational invariance allows one to partition a conserved quantity into additive components, which often provides insight into underlying process. 

If we can understand these simple shift, stretch, and rotational invariances, we will understand much about the intrinsic structure of pattern. An explanation of natural pattern often means an explanation of how particular processes lead to particular forms of invariance.

\section{Background}

This section introduces basic concepts and notation. I emphasize qualitative aspects rather than detailed mathematics. The final section of this article provides historical background and alternative perspectives.

\subsection{Probability increments}

Define $q(z)\equiv\qz$ such that the probability associated with $z$ is $\qz\GD\Gpz$. This probability is the area of a rectangle with height $\qz$ and incremental width $\GD\Gpz$. 

The total probability is constrained to be one, as the sum of the rectangular areas over all values of $z$, which is $\sum\qz\GD\Gpz=1$. When the $z$ values are discrete quantities or qualitative labels for events, then the incremental measure is sometimes set to one everywhere, $\GD\Gpz\equiv 1$, with changes in the measure $\GD\Gpz$ made implicitly by adjusting $\qz$. The conservation of probability becomes $\sum\qz=1$.

If a quantitative scale $z$ has values that are close together, then the incremental widths are small, $\GD\Gpz\rightarrow\dGpz$, and the distribution becomes essentially continuous in the limit. The probability around each $z$ value is $\qz\,\dGpz$. Writing the limiting sum as a integral over $z$, the conservation of total probability is $\int\qz\,\dGpz=1$. 

The increments may be constant-sized steps $\dGpz=\dz$ on the $z$ scale, with probabilities $\qz\dGpz=\qz\dz$ in each increment. One may transform $z$ in ways that alter the probability expression, $\qz$, or the incremental widths, $\dGpz$, and study how those changes alter or leave invariant properties associated with the total probability, $\qz\dGpz$.

\subsection{Parametric scaling relations}

A probability pattern, $\qz\dGpz$, may be considered as a parametric description of two scaling relations, $\qz$ and $\Gpz$, with respect to the parameter $z$. Geometrically, $\qz\dGpz$ is a rectangular area defined by the parametric height, $\qz$, with respect to the parameter, $z$, and the parametric width, $\dGpz$, with respect to the parameter, $z$.  

We may think of $z$ as a parameter that defines a curve along the path $\lrp{\Gpz,\qz}$, relating a scaled input measure, $\Gpz$, to a scaled output probability, $\qz$. The followings sections describe how different invariances constrain these scaling relations.

\section{Shift invariance and the exponential form}

I show that shift invariance and the conservation of total probability lead to the exponential form of probability distributions in \Eq{boltz}. Thus, we may consider the main conclusions of statistical mechanics and maximum entropy as secondary consequences that follow from the primacy of shift invariance and conserved total probability.

\subsection{Conserved total probability}\label{conservedProb}

This section relates shift invariance to the conservation of total probability. Begin by expressing probability in terms of a transformed scale, $z\mapsto\trz$, such that $\qz = k_0f(\trz)$ and
\begin{equation*}
  \int\qz\dGpz = \int k_0f(\trz)\dGpz=1.
\end{equation*}
The term $k_0$ is independent of $z$ and adjusts to satisfy the conservation of total probability.

If we assume that the functional form $f$ is invariant to a shift of the transformed scale by a constant, $a$, then by the conservation of total probability
\begin{equation}\label{eq:ashiftTotal}
  \int k_0f(\trz)\dGpz = \int k_a f(\trz+a)\dGpz = 1.
\end{equation}
The proportionality constant, $k_a$, is independent of $z$ and changes with the magnitude of the shift, $a$, in order to satisfy the constraint on total probability. 

Probability expressions, $q(z)\equiv\qz$, are generally not shift invariant with respect to the scale, $z$. However, if our transformed scale, $z\mapsto\trz$ is such that we can write \Eq{ashiftTotal} for any magnitude of shift, $a$, solely by adjusting the constant, $k_a$, then the fact that the conservation of total probability sets the adjustment for $k_a$ means that the condition for $\trz$ to be a shift invariant canonical scale for probability is
\begin{equation}\label{eq:ashift}
  \qz = k_0f(\trz)=k_a f(\trz+a),
\end{equation}
which holds over the entire domain of $z$. 

The key point here is that $k_a$ is an adjustable parameter, independent of $z$, that is set by the conservation of total probability. Thus, the conservation of total probability means that we are only required to consider shift invariance in relation to the proportionality constant $k_a$ that changes with the magnitude of the shift, $a$, independently of the value of $z$. Appendix A provides additional detail about the conservation of total probability and the shift-invariant exponential form.

\subsection{Shift-invariant canonical coordinates}

This section shows the equivalence between shift invariance and the exponential form for probability distributions. 

Let $x\equiv\trz$, so that we can write the shift invariance of $f$ in \Eq{ashift} as
\begin{equation*}
  f(x+a)=\Ga_a f(x).
\end{equation*}
By the conservation of total probability, $\Ga_a$ depends only on $a$ and is independent of $x$. 

If the invariance holds for any shift, $a$, then it must hold for an infinitesimal shift, $a=\Ge$. By Taylor series, we can write
\begin{equation*}
  f(x+\,\Ge) = f(x)+\Ge f'(x) = \Ga_\Ge f(x).
\end{equation*}
Because $\Ge$ is small and independent of $x$, and $\Ga_0=1$, we can write $\Ga_\Ge=1-\Gl\Ge$ for a constant $\Gl$. Then the previous equation becomes 
\begin{equation*}
  f'(x)=-\Gl f(x).
\end{equation*}
This differential equation has the solution
\begin{equation*}
  f(x) = \hat{k}e^{-\Gl x},
\end{equation*}
in which $\hat{k}$ may be determined by an additional constraint. Using this general property for shift invariant $f$ in \Eq{ashift}, we obtain the classical exponential-Boltzmann form for probability distributions in \Eq{boltz} as
\begin{equation}\label{eq:trz}
  \qz = ke^{-\Gl\trz}
\end{equation}
with respect to the canonical scale, $\trz$. Thus, expressing observations on the canonical shift-invariant scale, $z\mapsto\trz$, leads to the classical exponential form. If one accepts the primacy of invariance, the ``energy,'' $\E$, of the Boltzmann form in \Eq{boltz} arises as a particular interpretation of the generalized shift-invariant canonical coordinates, $\trz$.

\subsection{Entropy as a consequence of shift invariance}

The transformation to obtain the shift-invariant coordinate $\trz$ follows from \Eq{trz} as
\begin{equation*}
  -\log\qz = \Gl\trz - \log k.  
\end{equation*}
This logarithmic expression of probability leads to various classical definitions of entropy and information \autocite{cover91elements,feynman98statistical}. Here, the linear relation between the logarithmic scale and the canonical scale follows from the shift invariance of probability with respect to the canonical scale, $\trz$, and the conservation of total probability. 

I interpret shift invariance and the conservation of total probability as primary aspects of probability patterns. Entropy and information interpretations follow as secondary consequences. 

One can of course derive shift invariance from physical or information theory perspectives. My only point is that such extrinsic concepts are unnecessary. One can begin directly with shift invariance and the conservation of total probability.

\subsection{Example: the gamma distribution}

Many commonly observed patterns follow the gamma probability distribution, which may be written as
\begin{equation*}
  \qz = kz^{\Ga\Gl}e^{-\Gl z}.
\end{equation*}
This distribution is not shift invariant with respect to $z$, because $z\mapsto a+z$ alters the pattern 
\begin{equation*}
  \qz = kz^{\Ga\Gl}e^{-\Gl z} \ne k_a(a+z)^{\Ga\Gl}e^{-\Gl (a+z)}.
\end{equation*}
There is no value of $k_a$ for which this expression holds for all $z$.

If we write the distribution in canonical form
\begin{equation}\label{eq:gamma}
  \qz=ke^{-\Gl\trz}=ke^{-\Gl\lrp{z-\Ga\log z}},
\end{equation}
then the distribution becomes shift invariant on the canonical scale, $\trz=z-\Ga\log z$, because $\trz\mapsto a+\trz$ yields
\begin{equation*}
  \qz=ke^{-\Gl(a+\trz)}=k_ae^{-\Gl\trz},
\end{equation*}
with $k_a=ke^{-\Gl a}$. Thus, a shift by $a$ leaves the pattern unchanged apart from an adjustment to the constant of proportionality that is set by the conservation of total probability. 

The canonical scale, $\trz=z-\Ga\log z$, is log-linear. It is purely logarithmic for small $z$, purely linear for large $z$, and transitions between the log and linear domains through a region determined by the parameter $\Ga$. 

The interpretation of process in relation to pattern almost always reduces to understanding the nature of invariance. In this case, shift invariance associates with log-linear scaling. To understand the gamma pattern, one must understand how process creates a log-linear scaling relation that is shift invariant with respect to probability pattern\autocite{frank10measurement,frank11a-simple,frank14how-to-read}.

\section{Stretch invariance and average values}

\subsection{Conserved average values}

Stretch invariance means that multiplying the canonical scale by a constant, $\trz\mapsto b\trz$, does not change probability pattern. This condition for stretch invariance associates with the invariance of the average value. 

To begin, note that for the incremental measure $\dGpz=\dtrz$, the constant in \Eq{trz} to satisfy the conservation of total probability is $k=\Gl$, because
\begin{equation*}
  \int_0^\infty\Gl e^{-\Gl\trz}\,\dtrz = 1,
\end{equation*}
when integrating over $\trz$.

Next, define $\ave{X}{\Gp}$ as the average value of $X$ with respect to the incremental measure $\dGpz$. 
Then the average of $\Gl\trz$ with respect to $\dtrz$ is
\begin{equation}\label{eq:aveT}
  \Gl\ave{\tr}{\tr} = \int \Gl^2\,\trz\, e^{-\Gl\trz}\,\dtrz = 1.
\end{equation}
The parameter $\Gl$ must satisfy the equality. This invariance of $\Gl\ave{\tr}{\tr}$ implies that any stretch transformation $\trz\mapsto b\trz$ will be canceled by $\Gl\mapsto\Gl/b$. See Appendix A for further details.

We may consider stretch invariance as a primary attribute that leads to the invariance of the average value, $\Gl\ave{\tr}{\tr}$. Or we may consider invariance of the average value as a primary attribute that leads to stretch invariance.

\subsection{Alternative measures}

Stretch invariance holds with respect to alternative measures, $\dGpz\ne\dtrz$. Note that for $\qz$ in \Eq{trz}, the conservation of total probability fixes the value of $k$, because we must have
\begin{equation*}
  \int ke^{-\Gl\trz}\dGpz = 1.
\end{equation*}
The average value of $\Gl\trz$ with respect to $\dGpz$ is
\begin{equation*}
  \int \Gl\trz ke^{-\Gl\trz}\dGpz = \Gl\ave{\tr}{\Gp}.
\end{equation*}
Here, we do not have any guaranteed value of $\Gl\ave{\tr}{\Gp}$, because it will vary with the choice of the measure $\dGpz$. If we assume that $\ave{\tr}{\Gp}$ is a conserved quantity, then $\Gl$ must be chosen to satisfy that constraint, and, from the fact that $\Gl\trz$ occurs as a pair, $\Gl\ave{\tr}{\Gp}$ is a conserved quantity. The conservation of $\Gl\ave{\tr}{\Gp}$ leads to stretch invariance, as in the prior section. Equivalently, stretch invariance leads to the conservation of the average value.

\subsection{Example: the gamma distribution}

The gamma distribution from the prior section provides an example. If we transform the base scale by a stretch factor, $z\mapsto bz$, then 
\begin{equation*}
  \qz=k_be^{-\Gl\lrp{bz-\Ga\log bz}}.
\end{equation*}
There is no altered value of $\Gl$ for which this expression leaves $\qz$ invariant over all $z$. By contrast, if we stretch with respect to the canonical scale, $\trz\mapsto b\trz$, in which $\trz=z-\Ga\log z$ for the gamma distribution, we obtain
\begin{equation*}
  \qz = ke^{-\Gl_b b\trz}=ke^{-\Gl\trz}
\end{equation*}
for $\Gl_b=\Gl/b$. Thus, if we assume that the distribution is stretch invariant with respect to $\dz$, then the average value $\Gl\ave{\tr}{z}=\Gl\angb{z-\Ga\log z}$ is a conserved quantity. Alternatively, if we assume that the average value
\begin{equation*}
  \Gl\ave{\tr}{z}=\Gl\angb{z-\Ga\log z} = \Gl\angb{z}-\Gl\Ga\angb{\log z}
\end{equation*}
is a conserved quantity, then stretch invariance of the canonical scale follows. 

In this example of the gamma distribution, conservation of the average value with respect to the canonical scale is associated with conservation of a linear combination of the arithmetic mean, $\angb{z}$, and the geometric mean, $\angb{\log z}$, with respect to the underlying values, $z$. In statistical theory, one would say that the arithmetic and geometric means are sufficient statistics for the gamma distribution.

\section{Consequences of shift and stretch invariance}

\subsection{Relation between alternative measures}

We can relate alternative measures to the canonical scale by $\dtrz=\tr'\dGpz$, in which $\tr' = |\dtrz/\dd\Gp|$ is the absolute value of the rate of change of the canonical scale with respect to the alternative scale. Starting with \Eq{aveT} and substituting $\dtrz=\tr'\dGpz$, we have
\begin{equation*}
  \Gl\ave{\tr\tr'}{\Gp} = \int\Gl^2\trz\tr' e^{-\Gl\trz}\dGpz = 1.
\end{equation*}
Thus, we recover a universally conserved quantity with respect to any valid alternative measure, $\dGpz$. 

\subsection{Entropy}

Entropy is defined as the average value of $-\log\qz$. From the canonical form of $\qz$ in \Eq{boltz}, we have 
\begin{equation}\label{eq:logq}
  -\log\qz =\Gl\trz - \log k.
\end{equation}
Average values depend on the incremental measure, $\dGpz$, so we may write entropy \autocite{wikipedia15partition} as
\begin{equation*}
  \ave{-\log\qz}{\Gp} = \ave{\Gl\trz - \log k}{\Gp} = \Gl\ave{\tr}{\Gp} - \log k_\Gp.
\end{equation*}
The value of $\log k_\Gp$ is set by the conservation of total probability, and $\Gl$ is set by stretch invariance. The value of $\ave{\tr}{\Gp}$ varies according to the measure $\dGpz$. Thus, the entropy is simply an expression of the average value of the canonical scale, $\tr$, with respect to some incremental measurement scale, $\Gp$, adjusted by a term for the conservation of total probability, $k$.

When $\Gp\equiv\tr$, then $k_\Gp=\Gl$, and we have the classic result for the exponential distribution
\begin{equation*}
  \ave{-\log\qz}{\tr} =  \Gl\ave{\tr}{\tr} - \log \Gl = 1-\log\Gl=\log e/\Gl,
\end{equation*}
in which the conserved value $\Gl\ave{\tr}{\tr}=1$ was given in \Eq{aveT} as a consequence of stretch invariance.

\subsection{Cumulative measure}

Shift and stretch invariance lead to an interesting relation between $-\log \qz$ and the scale at which probability accumulates. From \Eq{logq}, we have
\begin{equation*}
  -\frac{1}{\Gl}\,\dlog\qz = \dtrz=\tr'\dGpz.
\end{equation*}
Multiplying both sides by $\qz$, the accumulation of probability with each increment of the associated measure is
\begin{equation*}
  -\frac{1}{\Gl}\,\qz\,\dlog\qz = \qz\dtrz=\qz\tr'\dGpz.
\end{equation*}
The logarithmic form for the cumulative measure of probability simplifies to
\begin{equation*}
  -\frac{1}{\Gl}\,\qz\,\dlog\qz = -\frac{1}{\Gl}\,\dd\qz = \qz\dtrz.
\end{equation*}
This expression connects the probability weighting, $\qz$, for each incremental measure, to the rate at which probability accumulates in each increment, $\dd\qz=-\Gl\qz\dtrz$. This special relation follows from the expression for $\qz$ in \Eq{boltz}, arising from shift and stretch invariance and the consequent canonical exponential form.

\subsection{Affine invariance and the common scales}

Probability patterns are invariant to shift and stretch of the canonical scale, $\trz$. Thus, affine transformations $\trz\mapsto a+b\trz$ define a group of related canonical scales. In previous work, we showed that essentially all commonly observed probability patterns arise from a simple affine group of canonical scales \autocite{frank10measurement,frank11a-simple,frank14how-to-read}. This section briefly summarizes the concept of affine invariant canonical scales. Appendix B provides some examples.

A canonical scale $\tr(z)\equiv\tr$ is affine invariant to a transformation $G(z)$ if 
\begin{equation*}
  \tr\lrb{G(z)} = a + b\tr(z)
\end{equation*}
for some constants $a$ and $b$. We can abbreviate this notion of affine invariance as
\begin{equation}\label{eq:tcircg}
  \tr\circ G \sim \tr,
\end{equation}
in which ``$\sim$'' means affine invariance in the sense of equivalence for some constants $a$ and $b$. 

We can apply the transformation $G$ to both sides of \Eq{tcircg}, yielding the new invariance $\tr\circ G\circ G\sim\tr\circ G$. In general, we can apply the transformation $G$ repeatedly to each side any number of times, so that
\begin{equation*}
  \tr\circ G^{\,n}\sim\tr\circ G^{\,m}
\end{equation*}
for any nonnegative integers $n$ and $m$. Repeated application of $G$ generates a group of invariances---a symmetry group. Often, in practical application, the base invariance in \Eq{tcircg} does not hold, but asymptotic invariance
\begin{equation*}
  \tr\circ G^{(n+1)}\sim\tr\circ G^{\,n}
\end{equation*}
holds for large $n$. Asymptotic invariance is a key aspect of pattern \autocite{frank09the-common}.

\section{Rotational invariance and the Gaussian radial measure}

The following sections provide a derivation of the Gaussian form and some examples. This section highlights a few results before turning to the derivation.

Rotational invariance transforms the total probability $\qz\dtrz$ from the canonical exponential form into the canonical Gaussian form
\begin{equation}\label{eq:gauss}
  \Gl e^{-\Gl\trz}\dtrz \mapsto ve^{-\pi v^2\Rz^2}\,\dRz.
\end{equation}
This transformation follows from the substitution $\Gl\trz\mapsto \pi v^2\Rz^2$, in which the stretch invariant canonical scale, $\Gl\trz$, becomes the stretch invariant circular area, $\pi v^2\Rz^2$, with squared radius $v^2\Rz^2$. The new incremental scale, $v\dRz$, is the stretch invariant Gaussian radial measure. 

We can, without loss of generality, let $v=1$, and write $\A=\pi\Rz^2$ as the area of a circle. Thus the canonical Gaussian form
\begin{equation}\label{eq:gaussArea}
  \qz\dGpz = e^{-\A}\dRz
\end{equation}
describes the probability, $-\log\qz=\A$, in terms of the area of a circle, $\A$, and the incremental measurement scale, $\dGpz$, in terms of the radial increments, $\dRz$. 

Feynman \autocite{feynman98statistical} noted the relation between entropy, radial measure, and circular area. In my notation, that relation may be summarized as
\begin{equation*}
  \ave{-\log\qz}{\R} = \ave{\A}{\R}.
\end{equation*}
However, Feynman considered the circular expression of entropy as a consequence of the underlying notion of statistical mechanics. Thus, his derivation followed from an underlying canonical ensemble of particles

By contrast, my framework derives from primary underlying invariances. An underlying invariance of rotation leads to the natural Gaussian expression of circular scaling. To understand how rotational invariance leads to the Gaussian form, it is useful consider a second parametric input dimension, $\Gth$, that describes the angle of rotation \autocite{bryc95the-normal}. Invariance with respect to rotation means that the probability pattern that relates $q(z,\Gth)$ to $\Gp(z,\Gth)$ is invariant to the angle of rotation.

\subsection{Gaussian distribution}

I now show that rotational invariance transforms the canonical shift and stretch invariant exponential form into the Gaussian form, as in \Eq{gauss}. To begin, express the incremental measure in terms of the Gaussian radial measure as
\begin{equation*}
  \Gl\dtrz=\pi v^2\dRz^2=2\pi v^2\Rz\dRz,
\end{equation*}
from which the canonical exponential form $\qz\dtrz=\Gl e^{-\Gl\trz}\dtrz$ may be expressed in terms of the radial measure as
\begin{equation}\label{eq:Rmeasure}
  \Gl e^{-\Gl\trz}\dtrz = 2\pi v^2\Rz e^{-\pi v^2\Rz^2}\,\dRz.
\end{equation}

Rotational invariance means that for each radial increment, $v\dRz$, the total probability in that increment given in \Eq{Rmeasure} is spread uniformly over the circumference $2\pi v\Rz$ of the circle at radius $v\Rz$ from a central location. 

Uniformity over the circumference implies that we can define a unit of incremental length along the circumferential path with a fraction $1/2\pi v\Rz$ of the total probability in the circumferential shell of width $v\dRz$. Thus, the probability along an increment $v\dRz$ of a radial vector follows the Gaussian distribution
\begin{equation*}
  \lrp{1/2\pi v\Rz}\qz\dtrz = ve^{-\pi v^2\Rz^2}\,\dRz
\end{equation*}
invariantly of the angle of orientation of the radial vector. 

Here, the total probability of the original exponential form, $\qz\dtrz$, is spread evenly over the two-dimensional parameter space $\lrp{z,\Gth}$ that includes all rotational orientations. The Gaussian expression describes the distribution of probability along each radial vector, in which a vector intersects a constant-sized area of each circumferential shell independently of distance from the origin.

The Gaussian distribution varies over all positive and negative values, $\Rz\in\lrp{-\infty,\infty}$, corresponding to an initial exponential distribution in squared radii, $\Rz^2=\trz\in\lrp{0,\infty}$. We can think of radial vectors as taking positive or negative values according to their orientation in the upper or lower half planes. 

\subsection{Radial shift and stretch invariance}

The radial value, $\Rz$, describes distance from the central location. Thus, the average radial value is zero, $\ave{\R}{R}=0$, when evaluated over all positive and negative radial values. Shift invariance associates with no change in radial distance as the frame of reference shifts the location of the center of the circle to maintain constant radii. 

Stretch invariance associates with the conserved value of the average circular area
\begin{equation*}
  \Gl\ave{\tr}{\R} = \pi v^2\ave{\R^2}{\R} = \pi v^2\Gs^2 = \frac{1}{2},
\end{equation*}
in which the variance, $\Gs^2$, is traditionally defined as the average of the squared deviations from the central location. Here, we have squared radial deviations from the center of the circle averaged over the incremental radial measure, $\dRz$. 

When $\Gl=v^2=1$, we have $\Gs^2=1/2\pi$, and we obtain the elegant expression of the Gaussian as the relation between circular area and radial increments in \Eq{gaussArea}. This result corresponds to an average circular area of one, because $\ave{2\pi\R^2}{}=2\pi\Gs^2=1$. 

It is common to express the Gaussian in the standard normal form, with $\Gs^2=1$, which yields $v^2=1/2\pi$, and the associated probability expression obtained by substituting this value into \Eq{gauss}. 

\subsection{Transforming distributions to canonical Gaussian form}

Rotational invariance transforms the canonical exponential form into the Gaussian form, as in \Eq{gauss}. If we equate $\Rz=\sqrt{\trz}$ and $\Gl=\pi v^2$, we can write the Gaussian form as 
\begin{equation}\label{eq:gaussroot}
  \qz\dRz = \sqrt{\frac{\Gl}{\pi}}\, e^{-\Gl\trz}\dd\sqrt{\trz},
\end{equation}
in which 
\begin{equation*}
  \tilde{\Gs}^2 = \ave{\tr}{\sqrt{\tr}}
\end{equation*}
is a generalized notion of the variance. 

The expression in \Eq{gaussroot} may require a shift of $\trz$ so that $\trz\in\lrp{0,\infty}$, with associated radial values $\Rz=\pm\sqrt{\trz}$. The nature of the required shift is most easily shown by example.

\subsection{Example: the gamma distribution}

The gamma distribution may be expressed as $\qz\dGpz$ with respect to the parameter $z$ when we set $\trz=z-\Ga\log z$ and $\dGpz=\dz$, yielding
\begin{equation*}
  \qz\dz = ke^{-\Gl\lrp{z-\Ga\log z}}\dz,
\end{equation*}
for $z\ge0$. To transform this expression to the Gaussian radial scale, we must shift $\trz$ so that the corresponding value of $\Rz$ describes a monotonically increasing radial distance from a central location. 

For the gamma distribution, if we use the shift $\trz\mapsto\trz-\Ga=\lrp{z-\Ga\log z}-\Ga$ for $\Ga\ge0$, then the minimum of $\trz$ and the associated maximum of $\qz$ correspond to $\Rz=0$, which is what we need to transform into the Gaussian form. In particular, the parametric plot of the points $\lrp{\pm\Rz,\qz}$ with respect to the parameter $z\in\lrp{0,\infty}$ follows the Gaussian pattern. 

In addition, the parametric plot of the points $\lrp{\trz,\qz}$ follows the exponential-Boltzmann pattern. Thus we have a parametric description of the probability pattern $\qz$ in terms of three alternative scaling relations for the underlying parameter $z$: the measure $\dz$ corresponds to the value of $z$ itself and the gamma pattern, the measure $\dRz$ corresponds to the Gaussian radial measure, and the measure $\dtrz$ corresponds to the logarithmic scaling of $\qz$ and the exponential-Boltzmann pattern. Each measure expresses particular invariances of scale. 

\subsection{Example: the beta distribution}

A common form of the beta distribution is
\begin{equation*}
  \qz\dz = kz^{\Ga-1}(1-z)^{\Gb-1}\dz
\end{equation*}
for $z\in\lrp{0,1}$. We can express this distribution in canonical exponential form $ke^{-\Gl\trz}$ by the scaling relation
\begin{equation*}
  -\Gl\trz = (\Ga-1)\log z + (\Gb-1)\log(1-z),
\end{equation*}
with $\Gl>0$. For $\Ga$ and $\Gb$ both greater than one, this scaling defines a log-linear-log pattern\autocite{frank11a-simple}, in the sense that $-\Gl\trz$ scales logarithmically near the endpoints of zero and one, and transitions to a linear scaling interiorly near the minimum of $\trz$ at
\begin{equation}\label{eq:betaz}
  z^* = \frac{\Ga-1}{\Ga+\Gb-2}.
\end{equation}
When $0<\Ga<1$, the minimum (extremum) of $\trz$ is at $z^*=0$. For our purposes, it is useful to let $\Ga=\Gl$ for $\Gl>0$, and assume $\Gb>1$. 

Define $\tr^*$ as the value of $\trz$ evaluated at $z^*$. Thus $\tr^*$ is the minimum value of $\trz$, and $\trz$ increases monotonically from its minimum. If we shift $\trz$ by its minimum, $\trz\mapsto\trz-\tr^*$, and use the shifted value of $\trz$, we obtain the three standard forms of a distribution in terms of the parameter $z\in(0,1)$, as follows.

The measure $\dz$ and parametric plot $\lrp{z,\qz}$ is the standard beta distribution form, the measure $\dRz$ and parametric plot $\lrp{\pm\Rz,\qz}$ is the standard Gaussian form, and the measure $\dtrz$ and parametric plot $\lrp{\trz,\qz}$ is the standard exponential-Boltzmann form. 

\section{Rotational invariance and partitions}

The Gaussian radial measure often reveals the further underlying invariances that shape pattern. Those invariances appear from the natural way in which the radial measure can be partitioned into additive components.

\subsection{Overview}

Conserved quantities may arise from an underlying combination of processes. For example, we might know that a conserved quantity, $\R^2=x+y$, arises as the sum of two underlying processes with values $x$ and $y$. We do not know $x$ and $y$, only that their conserved sum is invariantly equal to $\R^2$. 

The partition of an invariant quantity into a sum may be interpreted as rotational invariance, because
\begin{equation*}
  \R^2 = x+y=\sqrt{x}^2 + \sqrt{y}^2
\end{equation*}
defines a circle with conserved radius $\R$ along the positive and negative values of the coordinates $\lrp{\sqrt{x},\sqrt{y}}$. That form of rotational invariance explains much of observed pattern, many of the classical results in probability and dynamics, and the expression of those results in the context of mechanics. 

The partition can be extended to a multidimensional sphere of radius $\R$ as
\begin{equation}\label{eq:rotate}
  \R^2= \sum\sqrt{x_i}^2.
\end{equation}
One can think of rotational invariance in two different ways. First, one may start with a variety of different dimensions, with no conservation in any particular dimension. However, the aggregate may satisfy a conserved total that imposes rotational invariance among the components.

Second, every conserved quantity can be partitioned into various additive components. That partition starts with a conserved quantity and then, by adding dimensions that satisfy the total conservation, one induces a higher dimensional rotational invariance. Thus, every conserved quantity associates with higher-dimensional rotational invariance. 

\subsection{Rotational invariance of conserved probability}

In the probability expression $\qz\dGpz$, suppose the incremental measure $\dGpz$ is constant, and we have a finite number of values of $z$ with positive probability. We may write the conserved total probability as $\sum_z\qz=1$. Then from \Eq{rotate}, we can write the conservation of total probability as a partition of $\R^2=1$ confined to the surface of a multidimensional sphere
\begin{equation*}
  \sum_z\sqrt{\qz}^2=1.
\end{equation*}
There is a natural square root spherical coordinate system, $\sqrt{\qz}$, in which to express conserved probability. Square roots of probabilities arise in a variety of fundamental expressions of physics, statistics, and probability theory \autocite{frieden04science,frank15dalemberts}.

\subsection{Partition of the canonical scale}

The canonical scale equals the square of the Gaussian radial scale, $\trz=\Rz^2$. Thus, we can write a two-dimensional partition from \Eq{rotate} as
\begin{equation*}
  \trz=\sqrt{x_1}^2+\sqrt{x_2}^2.
\end{equation*}
Define the two dimensions as
\begin{align*}
  \sqrt{x_1}&=w\equiv w(z,s)\\
  \sqrt{x_2}&=\wdot\equiv \wdot(z,s),
\end{align*}
yielding the partition for the canonical scale as
\begin{equation}\label{eq:tpart}
  \trz=w^2+\wdot^2.
\end{equation}
This expression takes the input parameter $z$ and partitions the resulting value of $\trz=\Rz^2$ into a circle of radius $\Rz$ along the path $\lrp{w,\wdot}$ traced by the parameter $s$. 

The radial distance, $\Rz$, and associated canonical scale value, $\trz=\Rz^2$, are invariant with respect to $s$. In general, for each dimension we add to a partition of $\trz$, we can create an additional invariance with respect to a new parameter. 

\subsection{Partition into location and rate}

A common partition separates the radius into dimensions of location and rate. Define $\wdot=\prt w/\prt s$ as the rate of change in the location $w$ with respect to the parameter $s$. Then we can use the notational equivalence $H_z\equiv\trz=\Rz^2$ to emphasize the relation to a classic expression in physics for a conserved Hamiltonian as
\begin{equation}\label{eq:ham}
  H_z=w^2+\wdot^2,
\end{equation}
in which this conserved square of the radial distance is partitioned into the sum of a squared location, $w^2$, and a squared rate of change in location, $\wdot^2$. The squared rate, or velocity, arises as a geometric consequence of the Pythagorean partitioning of a squared radial distance into squared component dimensions. Many extensions of this Hamiltonian interpretation can be found in standard textbooks of physics.

With the Hamiltonian notation, $H_z\equiv\trz$, our canonical exponential-Boltzmann distribution is
\begin{equation*}
  \qz\dd H_z = \Gl e^{-\Gl H_z}\dd H_z.
\end{equation*}
The value $H$ is often interpreted as energy, with $\dd H$ as the Gibbs measure. For the simple circular partition of \Eq{ham}, the total energy is often split into potential, $w^2$, and kinetic, $\wdot^2$, components.

In this article, I emphasize the underlying invariances and their geometric relations as fundamental. From my perspective, the interpretation of energy and its components are simply one way in which to describe the fundamental invariances. 

The Hamiltonian interpretation is, however, particularly useful. It leads to a natural expression of dynamics with respect to underlying invariance. For example, we can partition a probability pattern into its currently observable location and its rate of change
\begin{equation*}
  e^{-\Gl H_z} = e^{-\Gl w^2}e^{-\Gl \wdot^2}.
\end{equation*}
The first component, $w^2$, may be interpreted as the observable state of the probability pattern at a particular time. The second component, $\wdot^2$, may be interpreted as the rate of change in the probability pattern. Invariance applies to the combination of location and rate of change, rather than to either component alone. Thus, invariance does not imply equilibrium.

\section{Summary of invariances}

Probability patterns, $\qz$, express invariances of shift and stretch with respect to a canonical scale, $\trz$. Those invariances lead to an exponential form
\begin{equation*}
  \qz\dGpz = ke^{-\Gl\trz}\dGpz,
\end{equation*}
with respect to various incremental measures, $\dGpz$. This probability expression may be regarded parametrically with respect to $z$. The parametric view splits the probability pattern into two scaling relations, $\qz$ and $\Gpz$, with respect to $z$, forming the parametric curve defined by the points $\lrp{\Gpz,\qz}$. 

For the canonical scale, $\trz$, we may consider the sorts of transformations that leave the scale shift and stretch (affine) invariant, $\tr\circ G \sim \tr$, as in \Eq{tcircg}. Essentially all of the canonical scales of common probability patterns  \autocite{frank10measurement,frank11a-simple,frank14how-to-read} arise from the affine invariance of $\tr$ and a few simple types of underlying invariance with respect to $z$. 

For the incremental measure scale, $\dGpz$, four alternatives highlight different aspects of probability pattern and scale. 

The scale $\dz$ leads to the traditional expression of probability pattern, $\qz\dz$, which highlights the invariances that set the canonical scale, $\trz$. 

The scale $\dtrz$ leads to the universal exponential-Boltzmann form, $\qz\dtrz$, which highlights the fundamental shift and stretch invariances in relation to the conservation of total probability. 

This conservation of total probability may alternatively be described by a cumulative probability measure, $\dqz=-\Gl\qz\dtrz$. 

Finally, rotational invariance leads to the Gaussian radial measure, $\dRz$. That radial measure transforms many probability scalings, $\qz$, into Gaussian distributions, $\qz\dRz$. 

Invariances typically associate with conserved quantities \autocite{neuenschwander10emmy}. For example, the rotational invariance of the Gaussian radial measure is equivalent to the conservation of the average area circumscribed by the radial measure. That average circular area is proportional to the traditional definition of the variance. Thus, rotational invariance and conserved variance are equivalent in the Gaussian form.

The Gaussian radial measure often reveals the further underlying invariances that shape pattern. That insight follows from the natural way in which the radial measure can be partitioned into additive components.

\section{The primacy of invariance and symmetry}

\begin{quote}
It was Einstein who radically changed the \textit{way} people thought about nature, moving away from the mechanical viewpoint of the nineteenth century toward the elegant contemplation of the underlying symmetry principles of the laws of physics in the twentieth century (\textcite[p.~153]{lederman04symmetry}).
\end{quote}

The exponential-Boltzmann distribution in \Eq{boltz} provides the basis for statistical mechanics, Jaynesian maximum entropy, and my own invariance framework. These approaches derive the exponential form from different assumptions. The underlying assumptions determine how far one may extend the exponential-Boltzmann form toward explaining the variety of commonly observed patterns. 

I claim that one must begin solely with the fundamental invariances in order to develop a proper understanding of the full range of common patterns. By contrast, statistical mechanics and Jaynesian maximum entropy begin from particular assumptions that only partially reflect the deeper underlying invariances. 

\subsection{Statistical mechanics}

Statistical mechanics typically begins with an assumed, unseen ensemble of microscopic particles. Each particle is often regarded as identical in nature to the others. Statistical averages over the underlying microscopic ensemble lead to a macroscopic distribution of measurable quantities. The exponential-Boltzmann distribution is the basic equilibrium macroscopic probability pattern. 

In contrast with the mechanical perspective of statistical physics, my approach begins with fundamental underlying invariances (symmetries). Both approaches arrive at roughly the same intermediate point of the exponential-Boltzmann form. That canonical form expresses essentially the same invariances, no matter whether one begins with an underlying mechanical perspective or an underlying invariance perspective. 

From my point of view, the underlying mechanical perspective happens to be one particular way in which to uncover the basic invariances that shape pattern. But the mechanical perspective has limitations associated with the unnecessarily particular assumptions made about the underlying microscopic ensemble. 

For example, to derive the log-linear scaling pattern that characterizes the commonly observed gamma distribution in \Eq{gamma}, a mechanical perspective must make special assumptions about the interactions between the underlying microscopic particles.

Some may consider the demand for explicit mechanical assumptions about the underlying particles to be a benefit. But in practice, those explicit assumptions are almost certainly false, and instead simply serve as a method by which to point in the direction of the deeper underlying invariance that shapes the scaling relations and associated probability patterns. 

I prefer to start with the deeper abstract structure shaped by the key invariances. Then one may consider the variety of different particular mechanical assumptions that lead to the key invariances. Each set of particular assumptions that are consistent with the key invariances define a special case.

There have been many powerful extensions to statistical mechanics in recent years. Examples include generalized entropies based on assumptions about underlying particle mechanics \autocite{tsallis09introduction}, superstatistics as the average over heterogeneous microscopic sets \autocite{beck03superstatistics}, and invariance principles applied to the mechanical aspects of particle interactions \autocite{hanel11generalized}. 

My own invariance and scaling approach subsumes essentially all of those results in a simple and elegant way, and goes much further with regard to providing a systematic understanding of the commonly observed patterns \autocite{frank10measurement,frank11a-simple,frank14how-to-read}. However, it remains a matter of opinion whether an underlying mechanical framework based on an explicit microscopic ensemble is better or worse than a more abstract approach based purely on invariances. 

\subsection{Jaynesian maximum entropy}

Jaynes \autocite{jaynes57information,jaynes57informationII} replaced the old microscopic ensemble of particles and the associated mechanical entropy with a new information entropy. He showed that maximum entropy, in the sense of information rather particle mechanics, leads to the classic exponential-Boltzmann form.  A large literature extends the Jaynesian framework \autocite{presse13principles}. Axiomatic approaches transcend the original justifications based on intuitive notions of information \autocite{shore80axiomatic}.

Jaynes' exponential form has a kind of canonical scale, $\trz$. In Jaynes' approach, one sets the average value over the canonical scale to a fixed value, in our notation a fixed value of $\ave{\tr}{z}$. That conserved average value defines a constraint---an invariance---that determines the associated probability pattern \autocite{frank15maximum}. The Jaynesian algorithm is the maximization of entropy, subject to a constraint on the average value of some quantity, $\trz$. 

Jaynes struggled to go beyond the standard constraints of the mean or the variance. Those constraints arise from fixing the average values of $\trz=z$ or $\trz=z^2$, which lead to the associated exponential or Gaussian forms. Jaynes did discuss a variety of additional invariances \autocite{jaynes03probability} and associated probability patterns. But he never achieved any systematic understanding of the common invariances and the associated commonly observed patterns and their relations. 

I regarded Jaynes' transcendence of the particle-based microscopic ensemble as a strong move in the right direction. I followed that direction for several years \autocite{frank09the-common,frank10measurement,frank11a-simple,frank14how-to-read}. In my prior work, I developed the intrinsic affine invariance of the canonical scale, $\trz$, with respect to the exponential-Boltzmann distribution of maximum entropy. The recognition of that general affine invariance plus the variety of common invariances of scale \autocite{hand04measurement,luce08measurement} led to my systematic classification of the common probability patterns and their relationships \autocite{frank10measurement,frank11a-simple,frank14how-to-read}. 

In this article, I have taken the next step by doing away with the Jaynesian maximization of entropy. I replaced that maximization with the fundamental invariances of shift and stretch, from which I obtained the canonical exponential-Boltzmann form.

With the exponential-Boltzmann distribution derived from shift and stretch invariance rather than Jaynesian maximum entropy, I added my prior work on the general affine invariance of the canonical scale and the additional particular invariances that define the common scaling relations and probability patterns. We now have a complete system based purely on invariances.

\subsection{Conclusion}

Shift and stretch invariance set the exponential-Boltzmann form of probability patterns. Rotational invariance transforms the exponential pattern into the Gaussian pattern. These fundamental forms define the abstract structure of pattern with respect to a canonical scale. 

In a particular application, observable pattern arises by the scaling relation between the natural measurements of that application and the canonical scale. The particular scaling relation derives from the universal affine invariance of the canonical scale and from the additional invariances that arise in the particular application. 

Together, these invariances define the commonly observed scaling relations and associated probability patterns. The study of pattern often reduces to the study of how particular generative processes set the particular invariances that define scale. 

Diverse and seemingly unrelated generative processes may reduce to the same simple invariance, and thus to the same scaling relation and associated pattern. To test hypotheses about generative process and to understand the diversity of natural pattern, one must understand the central role of invariance. Although that message has been repeated many times, it has yet to be fully deciphered.

\section*{Acknowledgments}

\noindent National Science Foundation grant DEB--1251035 supports my research.  I began this work while on fellowship at the Wissenschaftskolleg zu Berlin.



\bigskip
\bibliography{main}


\appendix

\addcontentsline{toc}{section}{Appendix A: Technical issues and extensions}

\section*{Appendix A: Technical issues and extensions}

\subsection{Conserved total probability}

The relations between shift invariance and the conservation of total probability in Section \ref{conservedProb} form a core part of the article. Here, I clarify the particular goals, assumptions, and consequences. 

In Section \ref{conservedProb}, I assumed that the conservation of total probability and shift invariance hold. From those assumptions, \Eq{ashift} follows, and thus also the exponential-Boltzmann form of \Eq{trz}. 

I am not claiming that conservation of total probability by itself leads to shift invariance. Instead, my goal is to consider the consequences that follow from a primary assumption of shift invariance. 

The justification for a primary assumption of invariance remains an open problem at the foundation of much of modern physics. The opening quote from Weinberg expresses the key role of invariances and also the uncertainty about why invariances are fundamental. My only goal concerns the consequences that follow from the assumption of primary invariances.

\subsection{Conserved average values: \Eq{aveT}}

Below \Eq{aveT}, I stated that the average value $\Gl\ave{\tr}{\tr}=1$ remains unchanged after stretch transformation, $\trz\mapsto b\trz$. This section provides additional details. The problem begins with \Eq{aveT}, repeated here
\begin{equation*}
  \Gl\ave{\tr}{\tr} = \int \Gl^2\,\trz\, e^{-\Gl\trz}\,\dtrz = 1.
\end{equation*}
Make the substitution $\trz\mapsto b\trz$, which yields
\begin{equation*}
  \Gl b\ave{\tr}{\tr} = \int \Gl^2b^2\,\trz\, e^{-\Gl b\trz}\,\dtrz = 1,
\end{equation*}
noting that $\trz\mapsto b\trz$ implies $\dtrz\mapsto b\dtrz$, which explains the origin of the $b^2$ term on the right-hand side. Thus, \Eq{aveT} remains one under stretch transformation, implying that $\ave{\tr}{\tr}=1/\Gl b$.

\subsection{Primacy of invariance}

This article assumes the primacy of shift and stretch invariance. The article then develops the consequences of primary invariance. There are many other ways of understanding the fact that the foundational exponential-Boltzmann distribution expresses shift and stretch invariance, and the Gaussian distribution expresses rotational invariance. One can derive those invariances from other assumptions, rather than assume that they are primary.

Classical statistical mechanics derives shift and stretch invariance as consequences of the aggregate behavior of many particles. Jaynesian maximum entropy derives shift and stretch invariance as consequences of the tendency for entropy to increase plus the assumptions that total probability is conserved and that the average value of some measurement is conserved. In my notation, the conservation of $\angb{\Gl\trz}$ is equivalent to the assumption of stretch invariance. Often, this kind of assumption is similar to various conservation assumptions, such as the conservation of energy. 

Another way to derive invariance is by the classic limit theorems of probability. \textcite{gnedenko68limit} beautifully summarized a key aspect:
\begin{quote}
In fact, all epistemologic value of the theory of probability is based on this: that large-scale random phenomena in their collective action create strict, nonrandom regularity.
\end{quote}
The limit theorems typically derive from assumptions such as the summation of many independent random components, or in more complicated studies, the aggregation of partially correlated random components. From those assumptions, certain invariances may arise as consequences.

It may seem that the derivation of invariances from more concrete assumptions provides a better approach. But from a mathematical and perhaps ultimate point of view, invariance is often tautologically related to supposedly more concrete assumptions. For example, conservation of energy typically arises as an assumption in many profound physical theories. In those theories, one could chose to say that stretch invariance arises from conservation of energy or, equivalently, that conservation of energy arises from stretch invariance. It is not at all clear how we can know which is primary, because mathematically they are often effectively the same assumption.

My point of departure is the opening quote from Weinberg, who based his statement on the overwhelming success of 20th century physics. That success has partly (mostly?) been driven by studying the consequences that follow from assuming various primary invariances. The ultimate basis for those primary invariances remains unclear, but the profoundly successful consequences of proceeding in this way are very clear. These issues are very important. However, a proper discussion would require probing the basis of modern physics as well as many deep recent developments in mathematics, which is beyond my scope. I simply wanted to analyze what would follow from the assumption of a few simple primary invariances.

\subsection{Measurement theory}

Classical measurement theory develops a rational approach to derive and understand measurement scales \autocite{hand04measurement,luce08measurement}. Roughly speaking, a measurement scale is defined by the transformations that leave invariant the relevant relations of the measurement process. Different approaches develop that general notion of invariance in different ways or expand into broader aspects of pattern (e.g, \textcite{grenander96elements}). 

This article concerns probability patterns in relation to scale. The key is that probability patterns remain invariant to affine transformation, that is, to shift and stretch transformations. Thus different measurement scales lead to the same invariant probability pattern if they are affine similar. I discussed the role of affine similarity in several recent articles \autocite{frank10measurement,frank11a-simple,frank14how-to-read}. Here, I briefly highlight the main points.

Start with some notation. Let $\tr(z)\equiv\tr$ be a transformation of underlying observations $z$ that define a scale, $\tr$. Each scale $\tr$ has the property of being invariant to certain alterations of the underlying observations. Let a candidate alteration of the underlying observation be the generator, $\gr(z)\equiv\gr$. Invariance of the scale $\tr$ to the generator $\gr$ means that
\begin{equation*}
  \tr\left[\gr(z)\right] = \tr(z),
\end{equation*}
which we can write in simpler notation as
\begin{equation*}
  \tr\circ\gr = \tr.
\end{equation*}
Sometimes we do not require exact invariance, but only a kind of similarity. In the case of probability patterns, shift and stretch invariance mean that any two scales related by affine transformation $\tr = a + b\tr$ yield the same probability pattern. In other words, probability patterns are invariant to affine transformations of scale. Thus, with regard to the generator $\gr$, we only require that $\tr\circ\gr$ fall within a family of affine transformation of $\tr$. Thus, we write the conditions for two probability patterns to be invariant to the generator $\gr$ as
\begin{equation*}
  \tr\circ\gr = a + b\tr \sim \tr,
\end{equation*}
and thus the key invariance relation for probability patterns is affine similarity expressed as
\begin{equation*}
  \tr\circ\gr  \sim \tr,
\end{equation*}
which was presented in the text as \Eq{tcircg}. My prior publications fully developed this relation of affine similarity and its consequences for the variety of scales that define the commonly observed probability patterns  \autocite{frank10measurement,frank11a-simple,frank14how-to-read}. Appendix B briefly presents a few examples, including the linear-log scale.

\addcontentsline{toc}{section}{Appendix B: Invariance and the common canonical scales}

\section*{Appendix B: Invariance and the common canonical scales}

The variety of canonical scales may be understood by the variety of invariances that hold under different circumstances. I introduced the affine invariance of the canonical scale in \Eq{tcircg}. This section briefly summarizes further aspects of invariance and the common canonical scales. Prior publications provide more detail \autocite{frank10measurement,frank11a-simple,frank14how-to-read}. 

Invariance can be studied by partition of the transformation, $z\mapsto\trz$, into two steps, $z\mapsto w\mapsto \trz$. The first transformation expresses intrinsic invariances by the transformation $z\mapsto w(z)$, in which $w$ defines the new base scale consistent with the intrinsic invariances. 

The second transformation evaluates only the canonical shift and stretch invariances in relation to the base scale, $w\mapsto a +bw$. This affine transformation of the base scale can be written as $\tr(w)=a+bw$. We can define $\tr(w)\equiv\trz$, noting that $w$ is a function of $z$. 

\subsection{Rotational invariance of the base scale}

Rotational invariance is perhaps the most common base scale symmetry. In the simplest case, $w(z)=z^2$. If we write $x=z\cos\Gth$ and $y=z\sin\Gth$, then $x^2+y^2=z^2$, and the points $(x,y)$ trace a circle with a radius $z$ that is rotationally invariant to the angle $\Gth$. Many probability distributions arise from rotationally invariant base scales, which is why squared values are so common in probability patterns. For example, if $w=z^2$ and $\trz\equiv w$, then the canonical exponential form that follows from shift and stretch invariance of the rotationally invariant base scale is
\begin{equation*}
  \qz=ke^{-\Gl w}=ke^{-\Gl z^2},
\end{equation*}
which is the Gaussian distribution, as discussed in the text.

Note that the word \textit{rotation} captures an invariance that transcends a purely angular interpretation. Instead, we have component processes or measurements that satisfy an additive invariance constraint. For each final value, $z$, there exist a variety of underlying processes or outcomes that satisfy the invariance $\sum x_i^2=z^2$. 

The word \textit{rotation} simply refers to the diversity of underlying Pythagorean partitions that sum to an invariant Euclidean distance. The set of invariant partitions falls on the surface of a sphere. That spherical property leads to the expression of invariant additive partitions in terms of rotation.

\subsection{General form of base scale invariance}

The earlier sections established that the canonical scale of probability patterns is invariant to shift and stretch. Thus we may consider as equivalent any affine transformation of the base scale $w\mapsto a+bw$. 

We may describe additional invariances of $w$, such as rotational invariance, in the general form
\begin{equation}\label{eq:baseInv}
  w\circ G \sim w,
\end{equation}
in which $w\circ G\equiv w\lrb{G(z)}$. We read \Eq{baseInv} as: the base scale $w$ is invariant to transformation by $G$, such that $w\circ G = a+bw$ for some constants $a$ and $b$. The symbol ``$\sim$'' abbreviates the affine invariance of $w$. 

For example, we may express the rotational invariance of the prior section as
\begin{equation*}
  w(z,\Gth)=z^2(\cos^2\Gth+\sin^2\Gth)=z^2,
\end{equation*}
because $\cos^2\Gth+\sin^2\Gth=1$ for any value of $\Gth$. We can describe rotation by the transformation
\begin{equation*}
  G(z,\Gth)=(z,\Gth+\Ge),
\end{equation*}
so that the invariance expression is
\begin{equation*}
   w\circ G = w\lrb{G(z,\Gth)} = w(z,\Gth+\Ge)=z^2.
\end{equation*}
Thus, the base scale $w$ is affine invariant to the rotational transformation generator, $G$, as in \Eq{baseInv}. Although this form of rotational invariance seems trivial in this context, it turns out to be the basis for many classical results in probability, dynamics, and statistical mechanics.

\subsection{Example: linear-log invariance of the base scale}

The invariance expression of \Eq{baseInv} sets the conditions for base scale invariances. Although there are many possible base scales, a few dominate the commonly observed patterns \autocite{frank10measurement,frank11a-simple,frank14how-to-read}. In this article, I emphasize the principles of invariance rather than a full discussion of the various common scales.

Earlier, I discussed the log-linear scale associated with the gamma distribution. This section presents the inverse linear-log scale, which is
\begin{equation*}
  w(z) = \Ga\log(1+\Gb z).
\end{equation*}
When $\Gb z$ is small, $w$ is approximately $\Ga\Gb z$, which is linear in $z$. When $\Gb z$ is large, $w$ is approximately $\Ga\log(\Gb z)$, which is logarithmic in $z$. This linear-log scale is affine invariant to transformations
\begin{equation*}
  G(z) = \frac{(1+\Gb z)^\Ga-1}{\Gb},
\end{equation*}
because $w\circ G = \Ga w \sim w$. The transformation, $G$, is linear for small magnitudes of $z$ and power law for large magnitudes of $z$.

The linear-log base scale, $w$, yields the probability distribution 
\begin{equation*}
  \qz=ke^{-\Gl w}= k(1+\Gb z)^{-\Gg},
\end{equation*}
for $\Gg=\Gl\Ga$. This expression is the commonly observed Lomax or Pareto type II distribution, which is equivalent to an exponential-Boltzmann distribution for small $z$ and a power law distribution in the upper tail for large $z$.

We can combine base scales. For example, if we start with $w_1$, a rotationally invariant scale, $z\mapsto z^2$, and then transform those rotationally invariant values to a linear-log scale, $w_2$, we obtain $w_2\lrb{w_1(z)} = \Ga\log(1+\Gb z^2)$. This scale corresponds to the generalized Student's distribution
\begin{equation*}
  \qz=k(1+\Gb z^2)^{-\Gg}.
\end{equation*}
For small magnitudes of $z$, this distribution is linear in scale and Gaussian in shape. For large magnitudes of $z$, this distribution has power law tails. Thus, a rotationally invariant linear-log scale grades from Gaussian to power law as magnitude increases.

\subsection{The family of canonical scales}

The canonical scale, $\trz$, determines the associated probability pattern, $\qz=ke^{-\Gl\trz}$. What determines the canonical scale? The answer has two parts. 

First, each problem begins with a base scale, $w(z)\equiv w$. The base scale arises from the invariances that define the particular problem. Those invariances may come from observation or by assumption. The prior sections gave the examples of rotational invariance, associated with squared-value scaling, and linear to power-law invariance, associated with linear to log scaling. When the base scale lacks intrinsic invariance, we may write $w\equiv z$. Earlier publications provided examples of common base scales \autocite{frank10measurement,frank11a-simple,frank14how-to-read}.

Second, the canonical scale arises by transformation of the base scale, $\trz=\tr(w)$. The canonical scale must satisfy both the shift and stretch invariance requirements. If the base scale itself satisfies both invariances, then the base scale is the canonical scale, $\trz=w$. In particular, if the probability pattern remains invariant to affine transformations of the base scale $w\mapsto \Gd + \Gg w$, then the shift and stretch invariant distribution has the form
\begin{equation}\label{eq:qzw}
  \qz=ke^{-\Gl w}.
\end{equation}
Alternatively, $w$ may satisfy the shift invariance requirement, but fail the stretch invariance requirement \autocite{frank11a-simple,frank14how-to-read}. We therefore need to find a canonical transformation $\tr(w)$ that achieves affine invariance with respect to the underlying shift, $G(w)=\Gd+w$. The transformation
\begin{equation}\label{eq:trzw}
  \trz=\tr(w)=e^{\Gb w}
\end{equation}
changes a shift invariance of $w$ into a stretch invariance of $\trz$, because
\begin{equation*}
  \tr(\Gd+w)=e^{\Gb(\Gd+w)}=e^{\Gb\Gd}e^{\Gb w}=b\tr\sim\tr
\end{equation*}
for $b=e^{\Gb\Gd}$. We can write $\tr(\Gd+w)=\tr\circ G$, thus this expression shows that we have satisfied the affine invariance $\tr\circ G \sim \tr$ of \Eq{tcircg}.

Thus, shift invariance with respect to $w$ generates a family of scaling relations described by the parameter $\Gb$. The one parameter family of canonical scales in \Eq{trzw} expands the canonical exponential form for probability distributions to
\begin{equation}\label{eq:bw}
  \qz=ke^{-\Gl\trz}=ke^{-\Gl e^{\Gb w}}.
\end{equation}
The simpler form of \Eq{qzw} arises as a limiting case for $\Gb\rightarrow 0$. That limiting form corresponds to the case in which the base scale, $w$, is itself both shift and stretch invariant \autocite{frank11a-simple,frank14how-to-read}. Thus, we may consider the more familiar exponential form as falling within the expanded one parameter symmetry group of scaling relations in \Eq{trzw}. 

The expanded canonical form for probability patterns in \Eq{bw} and a few simple base scales, $w$, include essentially all of the commonly observed continuous probability patterns \autocite{frank11a-simple,frank14how-to-read}.

\subsection{Example: extreme values}

In some cases, it is useful to consider the probability pattern in terms of the canonical scale measure, $\dtrz=|\tr'|\dz$. Using $\trz=e^{\Gb w}$, distributions take on the form often found in the extreme value problems \autocite{frank11a-simple,frank14how-to-read}
\begin{equation*}
  \qz\dz = k w'e^{\Gb w-\Gl e^{\Gb w}}\dz,
\end{equation*}
in which $w'=|\dd w/\dz|$. For example, $w=z$ yields the Gumbel distribution, and $w=\log z$ yields the Fr{\'e}chet or Weibull form.

\subsection{Example: stretched exponential and L\'evy}

Suppose the base scale is logarithmic, $w(z)=\log z$. Then from \Eq{bw}, a candidate form for probability pattern is
\begin{equation}\label{eq:stretch}
  \qz=ke^{-\Gl z^\Gb}.
\end{equation}
This important distribution arises in various contexts \autocite{frank14how-to-read}, including the stretched exponential distribution and the Fourier domain spectral distribution that associates with the basic L\'evy distributions \autocite{frank09the-common}. 

In this case, the probability pattern is not shift and stretch invariant to changes in the value of $z$, because $z\mapsto \Gd+\Gg z$ changes the pattern. By contrast, if we start with the base scale $w=\log z$, then the probability pattern is shift and stretch invariant with respect to the canonical scale
\begin{equation*}
  \trz=e^{\Gb w} = z^\Gb,
\end{equation*}
because the affine transformation of the canonical scale, $z^\Gb\mapsto\Gd+\Gg z^\Gb$, does not alter the probability pattern in \Eq{stretch}, given that we adjust $k$ and $\Gl$ to satisfy the conservation of probability and the conservation of average value.

The way in which I presented these invariances may seem trivial. If we begin with \Eq{stretch}, then of course we have shift and stretch invariance with respect to $z^\Gb\mapsto\Gd+\Gg z^\Gb$. However, in practical application, we may begin with an observed pattern and then try to infer its structure. In that case, analysis of the observations would lead to the conclusion of shift and stretch invariance with respect to the canonical power law scaling, $z^\Gb$. 

Alternatively, we may begin with a theory that includes a complicated interaction of various dynamical processes. We may then ask what invariance property matches the likely outcome of those processes. The conclusion may be that, asymptotically, shift and stretch invariance hold with respect to $z^\Gb\mapsto\Gd+\Gg z^\Gb$, suggesting the power law form of the canonical scale. 

In general, the particular invariant canonical scale derives from observations or from assumptions about process. The theory here shows the ways in which basic required invariances strongly constrain the candidate canonical scales. Those generic constraints shape the commonly observed patterns independently of the special attributes of each problem.

\end{document}

%% file: altdistn-def.tex
\newcommand*{\qz}{q_z}
\newcommand*{\dq}{\dd q}
\newcommand*{\dqz}{\dq_z}

\newcommand*{\wdot}{\dot{w}}

\newcommand*{\dz}{\dd z}

\newcommand*{\dlog}{\dd\log}

\newcommand*{\Gpz}{\Gp_z}
\newcommand*{\dGpz}{\dd\Gpz}

\newcommand{\tr}{\mathrm{T}}
\newcommand{\gr}{\mathrm{G}}
\newcommand*{\R}{R}

\newcommand*{\trz}{\tr_z}

\newcommand*{\dtrz}{\dd\trz}

\newcommand*{\Rz}{\R_z}
\newcommand*{\dRz}{\dd\Rz}

\newcommand*{\E}{\mathcal{E}}
\newcommand*{\Ep}{\E^{\raisebox{2pt}{\scriptsize$\prime$}}}

\newcommand*{\A}{\GL}

\newcommand*{\ave}[2]{\angb{#1}_{#2}}

\newcommand*{\Ga}{\alpha}
\newcommand*{\Gb}{\beta}
\newcommand*{\Gd}{\delta}
\newcommand*{\GD}{\Delta}
\newcommand*{\Ge}{\epsilon}
\newcommand*{\Gg}{\gamma}

\newcommand*{\Gl}{\lambda}
\newcommand*{\GL}{\Lambda}

\newcommand*{\Gs}{\sigma}

\newcommand*{\Gth}{\theta}

\newcommand*{\Gp}{\psi}

\newcommand*{\angb}[1]{\left\langle#1\right\rangle}
\newcommand*{\lrp}[1]{\left(#1\right)}
\newcommand*{\lrb}[1]{\left[#1\right]}

\newcommand*{\dd}{\textrm{d}}

\newcommand*{\Eq}[1]{eqn~\ref{eq:#1}}

\newcommand*{\prt}{\partial}

\newcount\BoxNum \BoxNum 1\relax
\makeatletter
\newcommand*{\boxlabel}[1]{%
  \protected@write \@auxout {}{\string \newlabel {box:#1}{{\the\BoxNum}}{}}%
  \advance\BoxNum 1\relax}
\makeatother
